\documentclass[a4paper,12pt]{article}
\usepackage{amssymb,delarray,amsmath,stmaryrd}
\usepackage[dvips]{graphicx}

\input epsf.tex
\newdimen\xsize
\newdimen\oldbaselineskip
\newdimen\oldlineskiplimit
\def\restorelineskip{\baselineskip=\oldbaselineskip%
\lineskiplimit=\oldlineskiplimit}
\def\putm[#1][#2]#3{
\hbox{\vbox to 0pt{\parindent=0pt%
\vskip#2\xsize\hbox to0pt{\hskip#1\xsize $#3$\hss}\vss}}}%
\long\def\Line#1{\hbox to \hsize{#1}}
\def\putt[#1][#2]#3{
\vbox to 0pt{\noindent\hskip#1\xsize\lower#2\xsize%
\vtop{\restorelineskip#3}\vss}}

\makeatletter
\def\xbig[#1]#2{{\hbox{$\m@th\left#2\vbox to#1\xsize{}%
\right.\n@space$}}}
\makeatother
\def\xlar[#1]#2{%
\smash{\mathop{ \hbox to #1\xsize{\leftarrowfill}}\limits^{#2}}}
\def\xrar[#1]#2{%
\smash{\mathop{ \hbox to #1\xsize{\rightarrowfill}}\limits^{#2}}}
\def\xline[#1]{\hbox to #1\xsize{\leaders\hrule\hfill}}

\DeclareFontFamily{U}{rsf}{\skewchar\font'177}%
\DeclareFontShape{U}{rsf}{m}{n}{<-6>rsfs5<6-8>rsfs7<8->rsfs10}{}%
\DeclareFontShape{U}{rsf}{b}{n}{<-6>rsfs5<6-8>rsfs7<8->rsfs10}{}%
\DeclareMathAlphabet\RSFS{U}{rsf}{m}{n}
\SetMathAlphabet\RSFS{bold}{U}{rsf}{b}{n}
\def\mib#1{\boldsymbol{#1}}
\def\sf#1{{\mathsf{#1}}}

\def\slsf{\slshape \sffamily }

\def\msmall#1{\mathchoice{\hbox{\small$\displaystyle {#1}$}}{#1}{#1}{#1}}        

\let\3=\ss

\let\ss=\sss
\def\bb{{\mathbb B}}
\def\cc{{\mathbb C}}

\def\rr{{\mathbb R}}

\def\st{_{\mathsf{st}}}
\def\area{\sf{area}}

\def\cl{\sf{cl}}

\def\dim{\sf{dim}\,}
\def\dimc{\dim_\cc}

\def\ext{\sf{ext}}

\def\sb{\bf {S}}

\def\id{\sf{Id}}
\def\im{\sf{Im}\,}

\def\re{\sf{Re}\,}

\def\lim{\mathop{\sf{lim}}}

\def\eps{\varepsilon}

\def\<{\langle}\let\la=\<
\def\>{\rangle}\let\ra=\>
 \let\bs=\bss 
\def\comp{\Subset}
\def\d{\partial}
\def\dbar{{\barr\partial}}

\def\ddef{\mathrel{{=}\raise0.3pt\hbox{:}}}
\def\deff{\mathrel{\raise0.3pt\hbox{\rm:}{=}}}
\def\inv{^{-1}}

\def\fraction#1/#2{\mathchoice{{\msmall{ #1\over#2}}}%
{{ #1\over #2 }}{{#1/#2}}{{#1/#2}}}
\def\norm#1{\left\Vert{#1}\right\Vert}

\def\le{\leqslant}

\def\emptyset{\varnothing}
\def\scirc{\mathop{\mathchoice{\hbox{\small$\circ$}}{\hbox{\small$\circ$}}%
{{\scriptscriptstyle\circ}}{{\scriptscriptstyle\circ}}}}
\def\longpoints{\leaders\hbox to 0.5em{\hss.\hss}\hfill \hskip0pt}
\def\stateskip{\smallskip}
\def\state#1. {\stateskip\noindent{\bf#1. }} 
\def\statep#1. {\stateskip\noindent{\bf#1 }} 
\def\proof{\state Proof. \2}

\def\Chi{\raise 2pt\hbox{$\chi$}}

\def\ie{\hskip1pt plus1pt{\sl i.e.\/,\ \hskip1pt plus1pt}}

\def\sli{{\sl i)} } 
\def\slii{{\sl i$\!$i)} } 
\def\sliii{{\sl i$\!$i$\!$i)} }

\def\barr#1{\mskip1mu\overline{\mskip-1mu{#1}\mskip-1mu}\mskip1mu}
\def\Chi{\raise 2pt\hbox{$\chi$}}

\let\phI=\phi\let\phi=\varphi\let\varphi=\phI
\def\bfbeta{{\boldsymbol\beta}}%
\def\calc{{\cal C}}

\def\calm{{\cal M}}

\def\calo{{\cal O}}

\def\eps{\varepsilon}

\def\bs{\backslash}

\def\comp{\Subset}
\def\d{\partial}
\def\dbar{{\barr\partial}}
\def\1{{1\mkern-5mu{\rom l}}}

\def\ge{\geqslant}
\def\inv{^{-1}}

\def\fraction#1/#2{\mathchoice{{\msmall{ #1\over#2}}}%
{{ #1\over #2 }}{{#1/#2}}{{#1/#2}}}

\def\le{\leqslant}

\def\emptyset{\varnothing}
\newcommand{\2}{\thinspace}

\def\qed{\ \ \hfill\hbox to .1pt{}\hfill\hbox to .1pt{}\hfill $\square$\par}

\def\comment#1\endcomment{}

 \belowdisplayskip=\abovedisplayskip

\def\lineeqqno(#1){\hfill\llap{\vbox to 10pt%
{\vss\begin{align} \eqqno(#1)\end{align}\vss}}\vskip1pt}

\advance\headheight 1.2pt

\def\ShowwLLabel#1{}

\def\thechpt{\Roman{chpt}}
\def\newchapt[#1]#2{\newpage%
\refstepcounter{chpt}\setcounter{subsection}{0}%
\setcounter{thm}{0}\setcounter{defi}{0}%
\setcounter{rema}{0}\setcounter{exrc}{0}%
\renewcommand{\thesubsection}{\thechpt.\arabic{subsection}}%
\section*{\begin{center}\huge \bf Chapter \thechpt\\ 
#2 \end{center}}\label{#1}%
\ \smallskip%
\addcontentsline{toc}{part}{Chapter \thechpt. #2}%
\markboth{Chapter \thechpt}{#2}%
}
\def\newsect[#1]#2{\refstepcounter{section}\setcounter{equation}{0}%
\renewcommand{\thesubsection}{\arabic{section}.\arabic{subsection}}%
\section*{\arabic{section}.
#2}\vspace{-20pt}\label{#1}\vspace{20pt}%
\addcontentsline{toc}{chapter}{Section \arabic{section}. #2}%
\markboth{Section \arabic{section}}{#2}}

\def\newlect[#1]#2{\refstepcounter{section}%
\renewcommand{\thesubsection}{\arabic{section}.\arabic{subsection}}%
\section*{Lecture \arabic{section}\\
#2}\label{#1}%
\addcontentsline{toc}{chapter}{Lecture \arabic{section}. #2}%
\markboth{Lecture \arabic{section}}{#2}}

\def\newprg[#1]#2{\refstepcounter{subsection}%
\subsection*{{\mdseries\slshape \sffamily \thesubsection.\ #2}} \label{#1}%
\addcontentsline{toc}{section}{\thesubsection. #2}%
}
   
\def\newappx[#1]#2{%
\refstepcounter{appx}\setcounter{section}{0}%
\renewcommand{\thesubsection}{A\arabic{appx}.\arabic{subsection}}%
\section*{Appendix \arabic{appx}\\ #2}
\label{#1}%
\addcontentsline{toc}{chapter}{Appendix A\arabic{appx}. #2}
\markboth{Appendix A\arabic{appx}}{#2}
}

\newtheorem{thm}{Theorem}
   \def\newthm#1{\begin{thm}\label{#1}} 

\newtheorem{nnthm}{Theorem.} 
   \def\newnnthm#1{\begin{nnthm} \label{#1}}
\newtheorem{lem}{Lemma}[section]
   \def\newlemma#1{\begin{lem} \label{#1}}

\newtheorem{prop}{Proposition}[section]
   \def\newprop#1{\begin{prop}\label{#1}}

\newtheorem{corol}{Corollary}[section]
   \def\newcorol#1{\begin{corol} \label{#1}}

\newtheorem{defi}{Definition}
   \def\newdefi#1{\begin{defi} \label{#1}\rm }

\newtheorem{exmp}{Example}
   \def\newexmp#1{\begin{exmp} \label{#1}\rm }

\newtheorem{exrc}{Exercise}
   \def\newexrc#1{\begin{exrc} \label{#1}\rm }

\newtheorem{rema}{Remark}
   \def\newrema#1{\begin{rema} \label{#1}\rm }

\def\eqqno(#1){\label{(#1)}}
\def\eqqref(#1){(\ref{(#1)})}

\title{Schwarz Reflection Principle, Boundary Regularity and Compactness for $J$-Complex Curves}
\author{S. Ivashkovich, A. Sukhov}
\begin{document}
\maketitle
\newtheorem{quest}{Question}

\begin{abstract}
We establish the Schwarz Reflection Principle for $J$-complex discs
attached to a real analytic $J$-totally real submanifold of an
almost complex manifold with real analytic $J$. We also prove the
precise boundary regularity and derive the precise convergence in Gromov
compactness theorem in $\calc^{k,\alpha}$-classes
\footnote{Key-words: almost complex structure, totally real
manifold, holomorphic disc, reflection principle.}.
\end{abstract}

\newsect[sect1]{Introduction}

\newprg[prg1.1]{Reflection Principle}

Denote by  $\Delta$ the unit disc in $\cc$, by ${\bf S}$ - the unit
circle. Let $\beta\subset \sb $ be a non-empty open subarc of $\sb$.

\medskip\noindent{\bf Theorem 1} {\slsf (Reflection Principle).}
{\it Let $(X,J)$ be a real analytic almost complex manifold and $W$
a real analytic $J$-totally real submanifold of $X$. Let $u:\Delta
\to X$ be a $J$-holomorphic map continuous up to $\beta$ and such
that $u(\beta) \subset W$.  Then  $u$ extends to a neighborhood of
$\beta$ as a (real analytic) $J$-holomorphic map. }

\smallskip The case of integrable $J$ is due to H. A. Schwarz \cite{Sw}.
Indeed, one can find local holomorphic coordinates in a neighborhood
of $u(p)$ for a taken $p\in\beta$ such that $W=\rr^n$ in these
coordinates and now the Schwarz Reflection Principle applies. In our
case there is no such reflection, since a general almost complex
structure doesn't admits any local (anti)-holomorphic maps. But the
extension result still holds.

\smallskip One can put Theorem 1
into a more general form of Carath\'eodory, \cite{Ca}. For this
recall that the cluster set $\cl (u, \beta)$ of $u$ at $\beta$
consists of all limits $\lim_{k\to\infty} u(\zeta_k)$ for all
sequences $\{ \zeta_k\}\subset \Delta$ converging to $\beta$. In
\cite{CGS} it was proved that if the cluster set $\cl (u, \beta)$ of
a $J$-holomorphic map $u:\Delta\to X$ is compactly contained in a
totally real submanifold $W$  then $u$ smoothly extends to $\beta$.
Therefore we derive the following

\smallskip\noindent{\bf Corollary 1.} {\it In the conditions of the
Theorem 1 the assumption of continuity of $u$ up to $\beta$ and
$u(\beta)\subset W$ one can replace by the assumption that $\overline{u(\Delta)}$ is
compact and the cluster set $\cl (u, \beta)$ is  contained
in $W$. }

\newprg[prg1.2]{Boundary Regularity}

For the proof of our Reflection Principle we need to study not only
real analytic boundary values but also the smooth ones (with finite
smoothness). For our method to work we need the precise regularity and
a certain kind of uniqueness of smooth  $J$-complex discs attached to a
$J$-totally real submanifold. The result obtained is the following

\medskip\noindent{\bf Theorem 2} {\slsf (Boundary Regularity).}
{\it Let $u:(\Delta ,\beta)\to (X,W)$ be a $J$-holomorphic map of
class $L^{1,2}\cap \calc^0(\Delta\cup\beta)$, where $W$ is $J$-totally real.
Then:
\begin{itemize}
\item[(i)] for any integer $k\ge 0$ and real $0<\alpha <1$ if
$J\in \calc^{k,\alpha}$ and $W\in \calc^{k+1, \alpha}$ then $u$ is
of class $\calc^{k+1, \alpha}$ on $\Delta\cup\beta$;
\item[(ii)] for $k\ge 1$ the condition $u\in L^{1,2}\cap\calc^0(\Delta\cup\beta)$
and $u(\beta)\subset W$ can be replaced by the assumption that $\overline{u(\Delta)}$
is compact and the cluster set $\cl (u, \beta)$ is 
contained in $W$.
\end{itemize}
}
\begin{rema}\rm
If $J$ is of class $\calc^0$ and $W$ of $\calc^1$ then
$u\in\calc^{\alpha}$ up to $\beta$ for all $0<\alpha <1$. This was
proved in \cite{IS2}, Lemma 3.1.
\end{rema}

\smallskip
For integrable $J$ the result of Theorem 2 is due to E. Chirka
\cite{Ch}. For non-integrable $J$ weaker versions of this Theorem
were obtained in \cite{CGS, GS, MS}. Namely, the
$\calc^{k,\alpha}$-regularity of $u$ up to $\beta$ was achieved
there under the same assumptions. The precise inner regularity was obtained by J.-C. Sikorav  in \cite{Sk}.

\newprg[prg1.3]{Compactness Theorem}

The precise regularity of $J$-complex discs attached to a
$J$-totally real submanifolds of Theorem 2 allows also to get the
precise convergency in Gromov compactness theorem. We refer to
the Subsection 6.1 and to  \cite{IS1,IS2} for the relevant notions
and definitions.

\medskip\noindent{\bf Theorem 3} {\slsf (Compactness Theorem).}
{\it Let  a sequence $\{ J_n\} $ of almost complex structures of
class $\calc^{k,\alpha}$, $k\ge 0, 0<\alpha <1$, on a Riemannian
manifold $(X,h)$ converge on a compact subset $K\subset X$ in
$\calc^{k,\alpha}$-topology to an almost complex structure $J$. Let
a sequence ${\mib W}_n =\{ (W_i,f_{n,i})\}_{i=1}^m$ of $J_n$-totally
real immersed submanifolds of $X$ of class $\calc^{k+1,\alpha}$
converge in $\calc^{k+1,\alpha}$-topology  to a $J$-totally real
immersion ${\mib W}=\{(W_i, f_{i})\}_{i=1}^m$. Suppose that all
${\mib W}_n$ and ${\mib W}$ have only weak transverse
self-intersections. Let furthermore, $\{ (\bar C_n, u_n)\} $ be a sequence
of stable $J_n$-complex curves over $X$, parameterized by a fixed
oriented compact real surface with boundary  $\bar\Sigma = (\Sigma ,\d\Sigma)$,
such that:

\sli $u_n(C_n) \subset K$ and that there exists a constant $M$ such
that $\area [u_n (C_n)]\le M$ for all $n$;

\slii $(\bar C_n,u_n)$ satisfy the totally real boundary conditions
$({\mib W}_n, \bfbeta , {\mib u}^{(b)}_n)$ with
$u^{(b)}_{n,i}(\beta_i)\subset W_i$ for all $n$ and $i$.

\smallskip
Then there exits a subsequence  $\{( \barr C_{n_k}, u_{n_k} )\}$ of
$\{( \barr C_n, u_n )\}$ and para\-metri\-zations $\sigma_{n_k}:
\barr \Sigma \to \barr C_{n_k}$, such that $(C_{n_k}, u_{n_k},
\sigma_{n_k})$ converges in $\calc^{k+1,\alpha}$-topology up to
boundary to a stable $J$-complex curve $(\barr C, u, \sigma)$ over
$X$ and this $(\barr C, u)$ satisfies the totally real boundary
conditions $({\mib W},\bfbeta, {\mib u}^{(b)})$ with some
$\calc^{k+1,\alpha}$-continuous
maps $u^b_{k}: \beta_k \to W_i$.}

\medskip The novelty here with respect  \cite{MS, F} is that there is no loss in both of
regularity and convergency of curves up to the boundary. In  \cite{IS2} an analog of 
Theorem 3 was proved for the special case $k = \alpha = 0$.

\newprg[prg1.4]{Proofs}

The interior analyticity of $J$-holomorphic discs in analytic
almost complex manifolds follows from classical results on elliptic
regularity in the real analytic category, see, for instance,
\cite{BJS}. However the real analyticity up to the boundary  is not a consequence of 
the known results since we do not deal with a
boundary problem of the Dirichlet type. The direct application of
the reflection principle (in the form of Vekua, for example) also
leads to technical complications because of the non-linearity of the
Cauchy-Riemann operator on an almost complex manifold. So our
approach is different and is based on the reduction of the boundary
regularity to a non-linear Riemann-Hilbert type boundary-value problem.

\smallskip
This paper is organized in the following way.

\smallskip\noindent 1. In \S 3, using \cite{IS2} and \cite{GS} we prove 
Theorem 2. First  we establish the $\calc^{1,\alpha}$-regularity
of $u$ if $J\in \calc^{\alpha}$ and then, using a sort of "geometric bootstrap",
we obtain the $\calc^{k+1,\alpha}$-regularity of $u$ if $J\in \calc^{k, \alpha}$. 

\smallskip\noindent 2. In \S 4 we prove the solvability and
uniqueness of a Riemann-Hilbert type boundary-value problem in Sobolev classes - the principal new
tool of this paper. In \S 5 we adapt this method to the real analytic case. Then the uniqueness, both in Sobolev and in real
analytic categories gives the proof of the Reflection Principle of
Theorem 1. The novelty in  \S 4 and 5, is a non-standard
choice of the smoothness classes and an application  of the Riemann-Hilbert boundary-value problem in the resolving
of boundary regularity questions.

\smallskip\noindent 3. In \S 6 we recall necessary notions and prove the Compactmess Theorem.

\smallskip\noindent 4. We end up with the formulation of open questions
in \S 7.

\medskip We would like to express our gratitude to
J.-F. Barraud who turned our attention to the question of
extendability of $J$-holomorphic maps through totally real
submanifolds in real analytic category.

\newsect[sect2]{Preliminaries}

In what follows $(X,J)$ will denote a pair which consists of a
real analytic manifold $X$ and an almost complex structure $J$ on
it. Note that by the well-known theorem of Whitney any smooth manifold
can be endowed with a compatible atlas with real analytic transition mappings, \ie the condition
of real analyticity of $X$ doesn't leed to any loss of generality.
The regularity of $J$ will be specified in each statement.  By
$J_{st}=i\id$ denote the standard complex structure of $\cc^n$ (as well as
of $\cc$ and of $\Delta$). Everywhere throughout the paper we use
the notation $L^{k,p}$ for Sobolev space of functions with
generalized partial derivatives of class $L^p$ up to the order $k$.
A $\calc^1$-map (or $\calc^0\cap L^{1,1}_{\sf loc}$-map)
$u:\Delta\to X$ is called {\slsf $J$-holomorphic} if for every (or
almost every) $\zeta\in\Delta$
\begin{eqnarray}
\eqqno(J-holomorphy1)
du(\zeta) \circ J_{st} = J(u(\zeta)) \circ
du(\zeta)
\end{eqnarray}
as mappings of tangent spaces $T_{\zeta}\Delta\to T_{u(\zeta)}X$.
The image $u(\Delta)$ is called then a $J$-{\slsf complex} disc. Every
almost complex manifold $(X,J)$ of complex dimension $n$ can be
locally viewed  as the unit ball $\bb$ in $\cc^n$ equipped with an
almost complex structure which is a small deformation of $J_{st}$.
To see this fix a point $p\in X$, choose a coordinate system such
that $p=0$, make an $\rr$-linear change of coordinates in order to
have $J(0)=J_{st}$ and rescale, i. e., consider $J(tz)$ for $t >0$
small enough. Then the equation \eqqref(J-holomorphy1) of
$J$-holomorphicity of a map $u:\Delta \to B$ can be written in local
coordinates $\zeta$ on $\Delta$ and $z$ on $\cc^n$ as the following
first order quasilinear system of partial differential equations
\begin{eqnarray}
\eqqno(J-holomorphy2)
u_{\overline\zeta} - A_J(u){\overline
u}_{\overline\zeta} = 0,
\end{eqnarray}
where $A_J(z)$ is the complex $n\times n$ matrix  of the operator
whose composite with complex conjugation is equal to the
endomorphism $ (J_{st} + J(z))^{-1}(J_{st} - J(z))$ (which is an
anti-linear operator with respect to the standard structure
$J_{st}$). Since $J(0) = J_{st}$, we have $A_J(0) = 0$. So in a
sufficiently small neighborhood of the origin the norm $\parallel A_J
\parallel_{L^\infty}$ is also small which implies the ellipticity
of the system \eqqref(J-holomorphy2).

We recall some classical integral transformations. Let  $\Omega$
be a relatively compact domain in $\cc$ bounded by a finite number
of transversally intersecting smooth curves.
Denote by $T^{CG}_\Omega$ the Cauchy-Green transform in $\Omega$:
\begin{equation}
\eqqno(CG-int)
\left( T^{CG}_\Omega h\right) (\zeta) = \frac{1}{2\pi i}\int\int_{\Omega}
\frac{h(\tau) d\tau \wedge d\overline \tau}{\tau - \zeta}.
\end{equation}
Denote also by
\begin{equation}
\eqqno(C-int)
\big(K_\Omega h\big)(\zeta) =\frac{1}{2\pi
i}\int_{\partial\Omega} \frac{h(\tau) d\tau }{\tau - \zeta}
\end{equation}
the Cauchy transform.

\begin{prop}
\label{CGoperator}
For every integer $k \geq 0$, real $0 < \alpha < 1$ and real $p>1$ the following holds:

\smallskip
(i) $T^{CG}_\Omega: \calc^{k, \alpha}(\bar\Omega) \longrightarrow
\calc^{k+1, \alpha}(\bar\Omega)$ (resp. $T^{CG}_\Omega: L^{k, p}(\Omega) \longrightarrow
L^{k+1, p}(\Omega)$) is a bounded linear operator and
$(T^{CG}_\Omega h)_{\overline\zeta} = h$ for any $h \in \calc^{k,
\alpha}(\Omega)$ (resp. for any $h \in L^{k, p}(\Omega)$).

\smallskip
(ii) Let $h$ be a bounded real analytic function on $\Omega$, then
$T^{CG}_\Omega h$ is real analytic on $\Omega$. If, furthermore, $\Omega = \Delta$,
$h$ is a sum of a converging series 
$h(\zeta,\overline\zeta) = \Sigma h_{k,l}\zeta^k\bar\zeta^l$ and  is continuous up to the boundary
$\d\Delta$, then for any $\zeta \in \Delta$ we have
\begin{eqnarray}
\label{analytic1}
T^{CG}_\Delta h(\zeta) = H(\zeta,\overline\zeta) -
\big(K_\Delta H\big)(\zeta),
\end{eqnarray}
where $H$ is a primitive of $h$ with respect to $\overline\zeta$. 
\end{prop}
The proof of (i) is contained, for instance, in \cite{Ve}, Theorems
1.32 and 1.37. For the statement (ii) about real analyticity see \cite{Ve}, p.26. 
Relation (\ref{analytic1}) is in fact nothing but the Cauchy-Green formula.  We shall also need the Schwarz integral
transform on $\Delta$:
\begin{equation}
\eqqno(SW-int)
\left( T^{SW} h\right) (\zeta) = \frac{1}{2\pi i} \int_{\partial\Delta} \frac{\tau + \zeta}{\tau - \zeta}\cdot
\frac{h(\tau)}{\tau}d\tau .
\end{equation}

Denote by $\calo^{k, \alpha}(\Omega)$ (resp. $\calo^{k, p}(\Omega)$) the Banach space of
holomorphic  maps $g: \Omega
\longrightarrow \cc^n$ of class $\calc^{k, \alpha}(\bar \Omega)$ (resp. $L^{k, p}(\bar \Omega)$).
This space is equipped with the norm $\parallel g
\parallel_{\calc^{k, \alpha}(\bar \Omega)}$ (resp. $\parallel g
\parallel_{L^{k, p}(\bar \Omega)}$).

\begin{prop}
\label{SWoperator} $T^{SW}:\calc^{k, \alpha}({\bf
S})\longrightarrow \calo^{k, \alpha}(\Delta)$ (resp. $T^{SW}:L^{k, p}({\bf
S})\longrightarrow \calo^{k, p}(\Delta)$) and $K_\Omega:
\calc^{k,\alpha}({\partial\Omega})\longrightarrow \calo^{k, \alpha}(\Omega)$ (resp. $K_\Omega:
L^{k,p}({\partial\Omega})\longrightarrow \calo^{k, p}(\Omega)$) are bounded
linear operators. For any real-va\-lued function $\psi \in
\calc^{k, \alpha}({\bf S})$  one has
\begin{equation}
\label{schwarz}
\re (T^{SW}\psi )|_{{\bf S}} = \psi
\end{equation}
and
\[
\im (T^{SW}\psi)(0) = 0.
\]
\end{prop}
For the proof  see, for instance, \cite{Ve}, Theorem 1.10. 

We shall
need to consider the space of traces of functions from
$L^{1,p}(\Delta)$ on the unit circle ${\bf S}$, see
\cite{Mo}. We say that a function $\varphi$ defined on ${\bf S}$ is
in the space $T^{1,p}(\bf S)$, $p
> 2$, if there exists a function $u \in L^{1,p}(\Delta)$ such that $u
\vert_{\bf S} = \varphi$. Let us point out that by the Sobolev imbedding
we have $T^{1,p}({\bf S}) \subset \calc^\alpha({\bf S})$ with $\alpha =
(p-2)/p$. On the other hand if $1/2 < \alpha <  1$ then
$\calc^{\alpha}({\bf S}) \subset T^{1,p}({\bf S})$ for every $p <
(1-\alpha)^{-1}$. Indeed, given $\varphi \in \calc^{\alpha}(\bf S)$
its Schwarz integral $u = T^{SW}\varphi$ is of class
$\calo^{\alpha}(\overline\Delta)$ and so by the classical theorem of
Hardy-Littlewood, see \cite{Mo}, one has
\[
\left\vert \frac{du(\zeta)}{d\zeta} \right\vert \leq \frac{C}{(1 -
\vert \zeta \vert)^{1 - \alpha}}.
\]
This implies that $\frac{du(\zeta)}{d\zeta} \in L^{p}(\Delta)$ for
$p(1-\alpha) < 1$. Now the Poisson integral
\begin{equation}
P\varphi:= \re T^{SW}(\re \varphi) + i\re T^{SW}(\im \varphi)
\label{Poisson}
\end{equation}
gives an extension of $\varphi$ of class $L^{1,p}(\Delta)$.

\begin{prop}
\label{SW} If $\varphi \in T^{1,p}({\bf S})$ then $T^{SW}\varphi \in
L^{1,p}(\Delta)$.
\end{prop}
\proof Let $u \in L^{1,p}(\Delta)$ be an extension of $\varphi$. By
the Cauchy-Green formula
\begin{equation}
u = K_\Delta \varphi + T^{CG}_{\Delta}u_{\overline\zeta}.
\label{Cauchy-Green}
\end{equation}
Hence the Cauchy type integral in the right hand (and therefore the Schwarz integral)
is of class $L^{1,p}(\Delta)$, which proves the proposition.

\smallskip\qed

Finally, we introduce the norm on the space $T^{1,p}({\bf S})$ by
setting
\[
\norm{ \varphi}_{T^{1,p}}:= \norm{P\varphi}_{L^{1,p}(\Delta)}
\]
where $P\varphi$ denotes the Poisson integral (\ref{Poisson}).
Obviously, $T^{1,p}({\bf S})$ is a Banach space. It follows from
(\ref{Poisson}) that the convergence in $\calc^{\alpha}({\bf S})$,
$\alpha > 1/2$ implies the convergence in $T^{1,p}({\bf S})$, \ie 
$\calc^{\alpha}({\bf S})$ is a closed subspace in $T^{1,p}({\bf S})$.
Furthermore, if a sequence of functions $\{u_n\}$ converges in
$L^{1,p}$, then their traces $\varphi_n:= u_n \vert_{\bf S}$
converge in $T^{1,p}({\bf S})$. This follows from the Cauchy-Green
representation (\ref{Cauchy-Green}) since the convergence of the
Cauchy type integral implies the convergence of the Schwarz and
Poisson integrals.

\begin{rema}\rm
In conclusion of this section we point out that the notion of the
trace space $T^{1,p}(\partial\Omega)$ on the boundary of a domain
$\Omega$ can be extended to a much larger class of simply connected
domains by means of the Riemann mapping theorem and the classical
theory of boundary properties of conformal mappings. For instance,
all above definitions and properties admits an immediate
generalization to the case where $\Delta$ is replaced by a simply
connected domain bounded by a finite number of real analytic arcs
with transversal intersections. The special case of the upper
semi-disc $\Delta^+$ will be important for our considerations. In
what follows we simply write $T^{1,p}$ in the case of the unit circle.
\end{rema}

\newsect[sect3]{Boundary Regularity in H\"oler Classes}

In this section we shall first use a version of a Reflection Principle
proposed in \cite{IS2} to prove the case $k=0$ of Theorem 2. Then using the "geometric bootstrap"
from \cite{GS} we obtain $\calc^{k+1,\alpha}$-regularity of complex discs in
$\calc^{k,\alpha}$-regular structures for all $k\ge 1$, thus proving the Theorem 2.

\newprg[prg3.1]{Reflection Principle-I}

In this subsection the structure $J$ is supposed to be of class $\calc^{\alpha}$ only. A $J$-totally 
real submanifold $W$ of $X$ will be supposed to have $\calc^{1,\alpha}$-regularity.
Similarly to the integrable case, a real submanifold $W$ of an almost complex manifold $(X,J)$ is 
called $J$-{\slsf totally real} if $T_pW\cap J\left( T_pW\right)  = \{ 0 \}$ at every point
$p$ of $W$. If $n$ is the complex dimension of $X$, then any totally real
submanifold of $X$ is locally contained in a totally real
submanifold of real dimension $n$. So in what follows we assume that
$W$ is $n$-dimensional. 

\smallskip First we make a suitable change of coordinates.

\begin{lem}
\label{alpha} 
One can find coordinates in a neighborhood $V$ of
$p\in W$ such that in these coordinates $V=\rr^{2n}$, $W=\rr^{n}$,
$J|_{\rr^n}=J\st$ and $J(x,y)-J\st = O(||y||^{\alpha})$.
\end{lem}
\proof After a change of coordinates of class $\calc^{1, \alpha}$ we
can suppose that in some neighborhood of $p=0$ our manifold $W$
coincides with $\rr^n$. Next we are looking for a $\calc^{1,
\alpha}$-diffeomorphism $\phi =(\phi_1,...\phi_{2n})$ in a
neighborhood of the origin such that

1) $\phi_j(x,0)=x_j$ for $j=1,...,n$;

2) $\phi_j(x,0)=0$ for $j=n+1,...,2n$;

3) $\frac{\d \phi}{\d y_j}(x,0) = J(x,0)\left(\frac{\d}{\d
x_j}\right)$ for $j=1,...,n$.

\noindent Such $\calc^{1, \alpha}$-diffeomorphism exists due to the
Trace theorem, see \cite{Tr}. In new coordinates given by $\phi$ we
shall clearly have $W=\rr^n$, $J|_{\rr^n}=J\st$ and $J(x,y)-J\st =
O(||y||^{\alpha})$ due to $\calc^{\alpha}$-regularity of $J$.

\qed

In what follows the disc $\Delta$ with an arc $\beta$ on its
boundary $\sb$ will be suitable for us to change by the upper
half-disc $\Delta^+=\{ \zeta : \re \zeta >0\}$ and the segment
$(-1,1)$. Let $u:(\Delta^+,\beta)\to (X,W)$ be a $J$-holomorphic map
of class $L^{1,p}$ up to $\beta = (-1,1)$ for some $p>2$.

\smallskip The following lemma will prove the case $k=0$ of
Theorem 2 and will be used in the proof of the same case of Theorem
3 in the last section.

\begin{lem}
\label{c1-alfa-reg} 
Let $J$ be of class $\calc^{\alpha}$ and $W$
be $J$-totally real of class $\calc^{1, \alpha}$. Let
$u:(\Delta^+,\beta)\to (X,W)$ be $J$-holomorphic of class
$\calc^0\cap L^{1,2}$ up to $\beta$. Then $u$ is of class $\calc^{1,
\alpha}$ up to $\beta$.
\end{lem}
\proof We can assume that $W=\rr^n$ and $J(x,y)-J\st =
O(||y||^{\alpha})$. On the trivial bundle  $\Delta^+\times
\rr^{2n}\to \Delta^+$ we consider the following linear complex
structure: $J_u(z)[\xi ]$ $= J(u(z))[\xi ]$ for $\xi \in\rr^{2n}$
and $z\in\Delta^+$. At this point we stress that $ J_u$ is defined
only on $\Delta^+\times\rr^{2n}$. Denote by $\tau$ the standard
conjugation in $\Delta\subset\cc$ as well as the standard
conjugation in $\rr^{2n}=\cc^n$. Now we extend $J_u$ to
$\Delta\times \rr^{2n}$ by setting
\begin{equation}
\tilde J_{u}(z)[\xi]= - \tau J_{u(\tau z)}[\tau \xi] \text{ for} z
\in \Delta \text{ and }\xi\in \rr^{2n}.
\end{equation}
We consider now $u$ as a section (over $\Delta^+$) of the trivial
bundle $E=\rr^{2n}\times\Delta\to\Delta$  and endow $E$ with the
complex structure $\tilde J_u$. Complex structure $\tilde J_u$
defines a $\dbar$-operator $\dbar_{\tilde J_u}w = \d_xw+\tilde
J_u\d_yw$ on $L^{1,p}$-sections of $E$ for all $1\le p<\infty$ (for
this only continuity of $\tilde J_u$ is needed). Remark that $u$ is
$\tilde J_u$-holomorphic on $\Delta^+$. By $F$ we denote the totally
real subbundle $\Delta\times \rr^n\to \Delta$ of $E$.

\begin{defi}
Define the ``extension by reflection'' operator $\ext
:L^1(\Delta^+,E)\to L^1(\Delta ,E)$ by setting
\begin{equation}
\label{extension} \ext (w)(z) = \tau w(\tau z)
\end{equation}
for $z\in \Delta^-$ and $w\in L^1(\Delta^+,E)$. We shall also write
$\tilde w$ for $\ext (w)$.
\end{defi}
Note that if $w$ is continuous up to $\beta$ and takes on $\beta$
values in the subbundle $F$ then $\ext (w)$ stays continuous. By the
reflection principle of Theorem 1.1 from \cite{IS2} we know that
$\ext : L^{1,p}(\Delta^+,E,F)\to L^{1,p}(\Delta ,E)$ is a continuous
operator for all $1\le p<\infty$ and that $\dbar_{\tilde J_u}\tilde
w = 0$ if $\dbar_{ J_u} w = 0$. Let $\tilde u$  be the extension of
$u$, in particular, $\tilde u$ is $\tilde J_u$-holomorphic on
$\Delta$ of class $L^{1,2}(\Delta)$. First a priori estimae (1.1)
from \cite{IS1} insures that $\tilde u\in L^{1,p}_{loc}$ for all
$p<\infty$. In particular $\tilde u\in \calc^{\gamma}$ for every
$\gamma <1$. Therefore $\tilde J_u$ is of class $\calc^{\delta}$,
where $\delta =\alpha\gamma$.

Set $v=\rho \tilde u$, where $\rho$ is a cut-off function equal to
$1$ in $\Delta_{\frac{1}{2}}$. Then from \eqqref(J-holomorphy2) we
get
\begin{equation}
v_{\bar\zeta} - A_J(\tilde u){\overline v}_{\overline\zeta} = g,
\end{equation}
where the function $g=\left[ \rho_{\bar\zeta} -
\rho_{\zeta}A_J(\tilde u)\right]\tilde u$ is of class
$\calc^{\delta}(\Delta)$. Observe that $v_{\bar\zeta} - A_J(\tilde
u){\overline v}_{\overline\zeta} = \left(
v-T^{CG}_{\Delta}A_J(\tilde u){\overline
v}_{\overline\zeta}\right)_{\overline\zeta}$. Elliptic regularity
implies that $v - T^{CG}_{\Delta}A_J(\tilde u){\overline
v}_{\overline\zeta}$ is of class $\calc^{1, \delta}$ and
invertibility of the operator $\id - T^{CG}_{\Delta}A_J(\tilde
u)\overline{\frac{\d}{\d z}}$ in $\calc^{1,\delta}$ gives that $v$
is in $\calc^{1,\delta}$.

One can repeat this step once more to get $\calc^{1,
\alpha}$-regularity of $v$ on $\Delta$ and therefore of $u$ up to
$\beta$.

\smallskip\qed

\newprg[prg3.2]{Cluster Sets on Totally Real Submanifolds}

Since the case $k=0$ of the Theorem 2 is already proved, we restrict
ourselves in the future with $k\ge 1$. Fix an almost complex
manifold $(X,J)$ with $J$ of class $\calc^{1, \alpha}$ and a
$J$-totally real submanifold $W$ of class $\calc^{2, \alpha}$. Let
$u: \Delta\to X$ be a bounded $J$-holomorphic map of the unit disc
into $X$. Suppose that $\cl (u,\beta)\comp W$, where $\beta $ is
some non-empty open subarc of the boundary.

\smallskip We use the Proposition 4.1 from \cite{CGS} and observe that $u$
is in Sobolev class $L^{1,p}$ up to $\beta$ for all $p<4$. In
particular $u$ is  $\calc^{\beta}$-regular up to $\beta$ with
$\beta=1-\frac{2}{p}$ (this means for all $\beta <\frac{1}{2}$).
Lemma 3.2 implies now the following

\begin{corol}
 Let $J\in \calc^{1,\alpha}$ and $W\in
\calc^{2,\alpha}$. If $u:(\Delta^+,\beta)\to (X,W)$ is a bounded
$J$-holomorphic map with $\cl (u,\beta)\comp W$ then
$u\in\calc^{1,\alpha}(\Delta^+\cup\beta)$.
\end{corol}

\medskip Let's stress here that $\calc^{1,\alpha}$ is not the
optimal regularity of $u$, it should be $\calc^{2,\alpha}$. This
will be achieved in the next subsection.

\newprg[prg3.3]{Geometric bootstrap and boundary $\calc^{k+1,\alpha}$-Regularity}

Having proved the $\calc^{1,\alpha}$-regularity of complex discs in
$\calc^{\alpha}$-regular structures we are going to use the "geometric bootstrap"
to obtain $\calc^{k+1,\alpha}$-regularity of complex discs in
$\calc^{k,\alpha}$-regular structures. 

\smallskip We  need to lift an
almost complex structure from the manifold $X$ to its tangent bundle $TX$.
In local coordinates the lift $J^c$ is defined by

\[
J^c=\left(
\begin{array}{ccc}
J_i^h & & 0\\
  & &    \\
t^a\partial_a J_i^h & & J_i^h
\end{array}
\right),
\]
where $t^a$ are coordinates in the tangent space. This lift is invariantly defined, see 
\cite{GS} for more details. After that we can use
the induction on $k$. Really, if $J\in\calc^{1,\alpha}$  then 
$J_n^c\in\calc^{\alpha}$. Further we  lift 
$J$-holomorphic map $u:\Delta\to X$  to  $J^c$-holomorphic map
$u^c:\Delta\to TX$. This lift is defined as

\begin{equation}
\eqqno(lift)
u^c(\zeta) = (u(\zeta),du(\zeta)(e_1)),
\end{equation}
where $e_1 = (1,0)$.

\smallskip Now remark that if $W$ is a $J$-totally real
submanifold in $X$ then $TW$ is a $J^c$-totally real
submanifold in $TX$. Really, let $v \in T(TW) \cap J^c(T(TW))$. If $v=(v_1,v_2)$ in the
trivialisation $T(TX) = TX\oplus TX$ then $v_1 \in TW \cap J(TW)$,
implying that $v_1=0$. Hence $v_2  \in TW \cap J(TW)$,
implying that $v_2=0$. Therefore $v=0$.

Further the lift $u^c:\Delta^+\to TX$ of a $J$-holomorphic map $u:\Delta^+\to X$ with boundary values in $W$ 
has its boundary values in $TW$ (this is clear from the formula \eqqref(lift)). 
Applying Lemma \ref{c1-alfa-reg} to $u^c$ and $TW$ we prove that the first derivative
of $u$ with respect to $\xi$
is of class $\mathcal \calc^{1,\alpha}$ on $\Delta^+ \cup \beta$. Here $\zeta = \xi + i\eta$.
The $J^c$-holomorphicity equation
\begin{equation}
\eqqno(J-holomorphy3)
\frac{\partial u^c}{\partial \eta} = J^c(u^c)
\frac{\partial u^c}{\partial \xi}
\end{equation}
implies that $\frac{\partial u^c}{\partial \eta}$ is also of class $\calc^{1,\alpha}$
on $\Delta^+ \cup \beta$. Therefore $u\in \calc^{2,\alpha}$ up to $\beta$. By induction we 
conclude that $u$ is of class $\calc^{k+1,\alpha}$ up
to $\beta$ if $J\in\calc^{k,\alpha}$ and $W$ is of class $\calc^{k+1,\alpha}$.

\smallskip Theorem 2 is proved.

\newsect[sect4]{Riemann-Hilbert Boundary-Value Problem}

In this Section we develop one of the main tools of this paper - a sort of 
a Riemann-Hilbert problem. This will be used in the proof of Theorem 1. Denote by ${\bf S}^+ = \{
e^{i\theta}: \theta \in ]0,\pi[ \}$ the upper semi-circle.  

\smallskip For the local conciderations of the present Section we suppose that 
$X=\cc^n$ and  $W = i\rr^n = \{ z = x +
iy: x = 0 \}$ and that $J\in\calc^{1, \alpha}$ is a small
deformation of $J_{st}$. Fix also a $J\st$-holomorphic map $u^0:
\Delta \longrightarrow \cc^n$ of class $L^{1, p}
(\Delta)$, $p > 2$, such that $u^0({\bf S}^+) \subset i\rr^n$ (so
that $u^0$ extends holomorphically to a neighborhood of ${\bf S}^+$
by the  classical Schwarz Reflection Principle). For $J$ close
enough to $J\st$ we will establish the existence and uniqueness of a $J$-holomorphic disc $u$ close enough to
$u^0$ satisfying the boundary condition $u({\bf S}^+)\subset
i\rr^n$.

Therefore for $J$ close enough to $J\st$ we  study the solutions of
\eqqref(J-holomorphy2) satisfying the boundary condition
\begin{eqnarray}
\label{realvanishing}
\re u|_{{\bf S}^+} = 0.
\end{eqnarray}
Denote by $T_0^{1, p}$  the Banach space of ($\rr^n$ -valued)
functions $\phi \in T^{1, p}$ vanishing on $\bf S^+$. This space is
equipped with the standard norm $\parallel \phi
\parallel_{T^{1, p}}$. Set now $\phi^0:= \re u^0
\vert {{\bf S}}$ and remark that $\phi^0\in T_0^{1, p}$ because $\re
u^0|_{\bf S^+}=0$. We replace the condition (\ref{realvanishing})
for the solutions of the partial differential equation
\eqqref(J-holomorphy2) by the condition
\begin{eqnarray}
\label{realequal}
\re u|_{\bf S} = \phi ,
\end{eqnarray}
where $\phi\in T_0^{1,p}$. Therefore we consider the  following
boundary-value problem
\begin{eqnarray}
\label{system1}
\begin{cases} u_{\overline\zeta} - A_J(u){\overline u}_{\overline\zeta} = 0, \cr
       \re u|_{\bf S} = \phi , \cr
       \im u(0) = a, \cr
\end{cases}
\end{eqnarray}
for the given initial data $\phi\in T^{1, p}_0$, $a\in \rr^n$.
\begin{lem}
\label{RH-smooth} If $J$ is close enough to $J\st$ in
$\calc^{1, \alpha}$-norm then there exists a neighborhood $U$ of $\phi^0$ in $T^{1,p}_0$,
a neighborhood ${ U}'$ of $a^0:=\im u^0(0)$ in $\rr^n$  and a
neighborhood $V$ of $u^0$  in $L^{1, p}(\Delta)$ such
that for each $\phi\in U$ and $a\in U'$ the boundary problem
(\ref{system1}) admits a unique solution $u \in V$.
\end{lem}
\proof   Consider the
operator
\begin{eqnarray*}
L_J: L^{1, p}(\Delta) \longrightarrow L^{p} (\Delta) \times
T^{1, p} \times \rr^n
\end{eqnarray*}
defined by
\begin{eqnarray*}
L_J: u \mapsto \left(
\begin{array}{cl}
 u_{\overline\zeta} - A_J(u){\overline u}_{\overline\zeta}\\
 \re u|_{\bf S}\\
 \im u(0)
\end{array}
\right) .
\end{eqnarray*}
$L_J$ smoothly depends on the parameter $J$. Denote by $\dot L_J(u)$
the  Fr\'echet derivative of $L_J$ at $u$. $\dot L_J$ is continuous on
the couple $(J,u)$ and at $J\st$-holomorphic $u^0$ the derivative
$\dot L_{J\st}(u^0)$ is particulary simple:

\begin{eqnarray*}
& &\dot L_{J\st}(u^0): L^{1, p}({\Delta}) \longrightarrow
L^{p}({\Delta})
\times T^{1, p}\times \rr^n\\
& &\dot L_{J\st}(u^0): \dot u \mapsto \left(
\begin{array}{cl}
 \dot u_{\overline\zeta}\\
 \re \dot u|_{\bf S}\\
 \im \dot u(0)
\end{array}
\right) .
\end{eqnarray*}
Let's see that $\dot L_{J\st}(u^0)$ is an isomorphism. Indeed, given
$h \in L^{p} ({\Delta})$, $\psi \in T^{1, p}$ and $a\in\rr^n$ then
the function
\[
\dot u = T^{CG}_{\Delta}h - i\im \left( T^{CG}_{\Delta}h(0)\right)
+ ia + T^{SW}(\psi -  \re (T^{CG}_{\Delta}h) \vert {\bf S})
\]
is of class $L^{1, p}({\Delta})$ and satisfies the equation
\begin{eqnarray*}
\dot L_{J\st}(u^0)(\dot u) = \left(
\begin{array}{cl}
 h\\
 \psi \\
 a
\end{array}
\right) .
\end{eqnarray*}
Uniqueness of $\dot u$ is obvious. Therefore by the Implicit
Function Theorem every  $L_J$ is a $\calc^1$-diffeomorphism of
neighborhoods of $u^0$ in $L^{1, p}(\Delta)$ and of $(0, \phi_0,
a_0)$ in $L^{p}({\Delta})\times T^{1, p}\times \rr^n$. Since $T^{1,
p}_0$ is a closed subspace of $T^{1, p}$ the Lemma \ref{RH-smooth}
follows.

\qed

Let us formulate a corresponding statement in H\"older classes. Let $k\ge 1$.
Fix a $J\st$-holomorphic map $u^0:
\Delta \longrightarrow \cc^n$ of class $\calc^{k+1, \alpha}
(\overline{\Delta})$ such that $u^0({\bf S}^+) \subset i\rr^n$. 
For every positive integer $k$ denote by $\calc_0^{k, \alpha}({\bf
S})$  the Banach space of ($\rr^n$ -valued) functions $\phi \in
\calc^{k, \alpha}({\bf S})$ vanishing on $\bf S^+$. This
space is equipped with the standard norm $\parallel \phi
\parallel_{\calc^{k, \alpha}(\sb)}$. Set now $\phi^0:= \re u^0
\vert {{\bf S}}$. $\phi^0\in \calc_0^{k+1, \alpha}({\bf S})$ because
$\re u^0|_{\bf S^+}=0$. Therefore we
consider the  boundary value problem (\ref{system1}) for the given initial data 
$\phi\in \calc^{k+1, \alpha}_0(\sb)$, $a\in \rr^n$.
\begin{lem}
\label{smoothness} Suppose $k\ge 1$. If $J$ is close enough to $J\st$ in
$\calc^{k, \alpha}$-norm then for every $1\le l\le k$:
\begin{itemize}
\item[(i)] there exists a neighborhood $U$ of $\phi^0$ in $\calc^{l+1,
\alpha}_0({\bf S})$,
a neighborhood ${ U}'$ of $a^0:=\im \phi^0(0)$ in $\rr^n$  and a
neighborhood $V$ of $u^0$  in $\calc^{l+1, \alpha}(\bar\Delta)$ such
that for each $\phi\in U$ and $a\in U'$ the boundary problem
(\ref{system1}) admits a unique solution $u \in V$;
\item[(ii)]  the unit disc $\Delta$ can be replaced in part (i) of the present
Lemma by any bounded simply connected domain $\Omega$ with
$\calc^\infty$ boundary and ${\bf S}^+$  can be replaced by any non-empty open
arc.
\end{itemize}
\end{lem}
\proof The part (ii) follows from (i) by the Riemann mapping theorem
and the classical theorems  on the boundary regularity of conformal
maps. For the proof of part consider the operator
\begin{eqnarray*}
L_J: \calc^{l+1, \alpha}(\bar\Delta) \longrightarrow \calc^{l, \alpha }
(\bar\Delta) \times \calc^{l+1, \alpha}({\bf S}) \times \rr^n ,
\end{eqnarray*}
defined as in the proof of Lemma \ref{RH-smooth}, but this time in other
smoothness clases. One can literally repeat the arguments used there also in
this case.

\smallskip\qed

\newsect[sect5]{Reflection Principle-II: Real Analytic Case}

Let us turn now to the proof of Theorem 1. We shall proceed in two
steps.

\newprg[prg5.1]{Small deformations of the standard structure} Here we
consider the case when $J$ is a small real analytic deformation of
$J_{st}$ and $W=i\rr^n$. First we introduce suitable Banach spaces
of real analytic functions using the complexification.

Denote by $\Delta^2 = \Delta \times \Delta$ the standard bidisc in
$\cc^2$. We define the space ${\calc^{1, \alpha}_\omega( \Delta)}$
consisting of functions $u$ (or  $\cc^n$-valued maps) of class $
{\calc^{1, \alpha}(\bar\Delta)}$ with the following properties:
\begin{itemize}
\item[(i)] $u$ is a sum of a power series $u(\zeta) = \sum_{kl}
u_{kl}\zeta^k\overline{\zeta}^l$ for $\zeta \in \Delta$.
\item[(ii)] The "polarization" $\hat u$ of $u$ defined by
 $\hat u(\zeta,\xi) = \sum_{kl} u_{kl}\zeta^k\xi^l$ is a function
 holomorphic
on $\Delta^2$ and of class ${\calc^{1, \alpha}(\overline
\Delta^2)}$. \item[(iii)] The mixed derivative $\frac{\partial^2
\hat u}{\partial \zeta \partial \xi}$ is of class ${\calc^{\alpha}( \overline\Delta^2)}$.
\end{itemize}

 We define the norm of $u$ as following:
\begin{eqnarray*}
\norm{u}_{\calc^{1,\alpha}_{\omega}(\Delta)} = \norm{\hat
u}_{\calc^{1,\alpha}(\overline \Delta^2)} + \norm{\frac{\partial^2
\hat u}{\partial \zeta \partial
\xi}}_{\calc^{\alpha}(\overline\Delta^2)}.
\end{eqnarray*}
Since $u$ is the restriction  of $\hat u$  onto the totally real
diagonal $\{ \xi =\bar\zeta\}$, the polarization $\hat u$ is
uniquely determined by $u$ and therefore $\calc^{1,\alpha}_\omega(\Delta)$
equipped with this norm is a Banach space.

\begin{rema}\rm
One has the following continuous inclusion
$\calo^{1,\alpha}(\Delta)\subset \calc^{1,\alpha}_{\omega}(\Delta)$:
for $u\in \calo^{1,\alpha}(\Delta)$ the corresponding $\hat u$ is
simply $\hat u(\zeta ,\xi )= u(\zeta)$. Really, for such $\hat u$
one has $\frac{\partial^2\hat u}{\partial \zeta \partial\xi}=0$.
\end{rema}

We denote by ${\calc}^{1, \alpha}_\omega({\d\Delta^+})$ the space
of real functions $\phi$ on $\d\Delta^+$ such that there exists a
function $v \in \calo^{1, \alpha }(\Delta)$ satisfying the
condition $\re v \vert_{\d\Delta^+} = \phi$. In particular such
function $\phi$ is real analytic on the interval $(-1,1)$. The
holomorphic function $v$ is unique up to an imaginary constant so
we always assume that $\im v(0) = 0$. We define the norm of $\phi$
as a $\calc^{1,\alpha}$ norm of the corresponding function $v$ on
$\bar\Delta$. Then ${\calc}^{1, \alpha}_\omega({\d\Delta^+})$
equipped with this norm, is  a Banach space.

Furthermore, denote by  $\calc^{1,\alpha}(\d\Delta^+)$ the space of
real continuous functions on $\d\Delta^+$ which are of class
$\calc^{1,\alpha}$ on the closed upper semi-circle and on the
interval $[-1,1]$. Finally we denote by
$\calc^{1,\alpha}_0(\d\Delta^+)$ the space of real functions of
class $\calc^{1,\alpha}(\d\Delta^+)$  vanishing on the interval
$[-1,1]$. The following statement is a consequence of the reflection
principle.
\begin{lem}
\label{reflectionlemma} 
For every function $\phi\in\calc^{1,\alpha}_0({\d\Delta^+})$ there exists $u\in \calo^{1,\alpha}(\Delta)$ 
such that $\re u|_{\d\Delta^+}=\phi$. In particular, the space $\calc^{1,\alpha}_0({\d\Delta^+})$ is a 
subspace of ${\calc}^{1,\alpha}_\omega({\d\Delta^+})$.
\end{lem}
\proof  Let $\phi$ be a function  of class
$\calc^{1,\alpha}_0(\d\Delta^+)$. Solving the Dirichlet problem
for $\phi$ in the upper semi-disc, we obtain a harmonic function
$h$ in $\Delta^+$ continuous on $\overline\Delta^+$ such that $h
\vert_{\d\Delta^+} = \phi$. Since $h$ vanishes on $[-1,1]$ it
extends harmonically on $\Delta$ by the classical reflection
principle for harmonic functions. Namely, its extension $h^*$ is
defined by $h^*(\zeta) = -h(\overline\zeta)$ for $\zeta$ in the
lower semi-disc $\Delta^-$. Thus we obtain a function $\tilde h$
harmonic on $\Delta$ and continuous on $\overline\Delta$. Since
the restriction $\phi$ of $h$ on the closed upper semi-circle is a
function of class $\calc^{1,\alpha}$, it follows easily  by the
definition of the reflection $h^*$ that  the restriction $\tilde
\phi:= \tilde h \vert_{\d\Delta}$ of $\tilde h$ on $\d\Delta$ is a
function of class $\calc^{1,\alpha}(\d\Delta)$. Then the Schwarz
integral $T^{SW}\tilde\phi$ gives by Proposition
\ref{SWoperator} a function of class $\calo^{1, \alpha
}(\Delta)$ whose real part coincides with $\tilde h$.

\qed

\begin{lem}
\label{image} If $u \in \calc_{\omega}^{1, \alpha }(\Delta)$ then
$\re u \vert_{\d\Delta^+} \in
\calc^{1,\alpha}_\omega({\d\Delta^+})$.
\end{lem}
\proof Let $\hat u(\zeta,\xi) = \sum u_{kl} \zeta^k\xi^l$ be the
polarization of $u$  holomorphic in the bidisc $\Delta^2$ (that is
$u(\zeta) = \hat u (\zeta,\overline \zeta)$ ). Then the function
$h(\zeta) = \hat u(\zeta,\zeta)$ is of class $\calo^{1, \alpha
}(\Delta)$ and $h \vert_{[-1,1]} = u \vert_{[-1,1]}$. Denote b
$\phi$ the restriction of $\re (u - h)$ to $\d\Delta^+$. Then $
\phi \in {\calc}^{1, \alpha}_0({\d\Delta^+})$ and by Lemma
\ref{reflectionlemma} there exists a function $v \in \calo^{1,
\alpha }(\Delta)$ such that $\re v \vert \d\Delta^+ = \phi$. Since
the function $ h + v$ is of class $\calo^{1, \alpha }(\Delta)$,
its real part gives the desired extension  of the function $\re u
\vert \d\Delta^+$.

\qed

We suppose everywhere below that our almost complex structure
$J$ (and therefore $A_J$ in the equation for $J$
holomorphic curves) is a real analytic matrix-valued function given by a
convergent power series $\sum a_{kl} z^k \overline{z}^l$ with the
radius of convergence big enough. The equation \eqqref(J-holomorphy2) on $\Delta$
can be rewritten in the form
\begin{equation}
\label{J-holomorphy3} (u - T^{CG}_{\Delta} A_J(u){\overline
u}_{\overline\zeta})_{\overline\zeta} = 0,
\end{equation}
where $T^{CG}_{\Delta}$ denotes the Cauchy - Green transform on
$\Delta$. Define the map
\[
\Phi_J : \calc^{1, \alpha}(\bar \Delta) \to \calc^{1, \alpha}(\bar
\Delta),
\]
as
\begin{equation}
\Phi_J : u \mapsto u - T^{CG}_{\Delta} A_J(u){\overline
u}_{\overline\zeta}.
\end{equation}

Equation (\ref{J-holomorphy3}) means that $u$ is $J$-holomorphic
if and only if  $\Phi_Ju$  is holomorphic with respect to
$J_{st}$. The following lemma explains the choice of smoothness
classes in this Section and is the principal step in the proof of Theorem 1.

\begin{lem}
\label{diffeo}
For $J$ close to $J\st$ the operator $\Phi_J$ establishes a
diffeomorphism of neighborhoods of zero in the space
$\calc^{1,\alpha}_\omega(\Delta)$.
\end{lem}
\proof  First we prove that $\Phi_J$ maps the space
$\calc^{1,\alpha}_\omega(\Delta)$ to itself. Given function $ u
\in \calc^{1,\alpha}_\omega(\Delta)$ denote the function
$A_J(u){\overline u}_{\overline\zeta}$ by $h$. We need to prove that
$T^{CG}_{\Delta}h$ belongs to $\calc^{1,\alpha}_\omega(\Delta)$. Consider
the polarization   $\hat h (\zeta,\xi) = \hat h(\zeta,\xi)$ of $h$.
By Proposition \ref{CGoperator} we have the representation
\begin{equation}
\label{lasteq} T^{CG}_\Delta  h(\zeta) = \hat
H(\zeta,\overline\zeta) - \frac{1}{2\pi i}\int_{\partial\Delta}
\frac{\hat H(\tau,\overline\tau) d\tau }{\tau - \zeta}.
\end{equation}
where
\begin{equation}
\label{primitive}
 \hat H(\zeta,\xi) = \int_{[0,\xi]} \hat
h(\zeta,\omega)d\omega
\end{equation}
is a primitive of $\hat h$ with respect to $\xi$.  Let's study the primitive
(\ref{primitive}) of $\hat h$ first. We point out that the function $\hat h$
is of class $\calc^{0,\alpha}(\overline\Delta^2)$. Furthermore, the condition
(iii) of the definition of the space $\calc^{1,\alpha}_\omega(\Delta)$ implies
that $\partial \hat h/\partial \zeta$ is of class $\calc^{0,\alpha}(\overline\Delta^2)$.
Now the derivation of  the integral (\ref{primitive}) with respect to $\zeta$
and $\xi$ gives that $\hat H$ satisfies conditions (i), (ii), (iii) of the definition
of the space $\calc^{1,\alpha}_\omega(\Delta)$.

By Proposition
\ref{SWoperator} the Cauchy integral in the right hand side of
(6.3) represents a function of class $\calo^{1,\alpha}(\Delta)$ and so also belongs
to the space $\calc^{1,\alpha}_\omega(\Delta)$.

Thus we obtain that $\Phi_J(u)$ belongs to
$\calc^{1,\alpha}_\omega(\Delta)$. Since the Fr\'echet derivative
of $\Phi_J$ with respect to $u$ at $u = 0$ and $J = J_{st}$ is the
identity map, the lemma follows from the inverse mapping theorem.

\smallskip\qed

\smallskip Hence $\Phi_J$ is a diffeomorphism
between neighborhoods of zero in the manifolds of $J$-holomorphic
and $J_{st}$-holomorphic maps of class $\calc^{1,
\alpha}_{\omega}(\Delta)$. In particular, $J$-holomorphic maps
form a Banach submanifold in $\calc^{1,\alpha}_\omega(\Delta)$ in a
neighborhood of zero. We denote this manifold as
$\calo^{1,\alpha}_{\omega ,J}(\Delta)$.

\begin{rema}\rm
Note that $\calo^{1,\alpha}_{\omega ,J\st}(\Delta)=\calo^{1,\alpha}(\Delta)$ and that 
$\calo^{1,\alpha}_{\omega ,J}(\Delta)=\Phi_J(\calo^{1,\alpha}(\Delta))$ .
\end{rema}
We use the notation
$\calo^{1,\alpha}_{\omega ,J,0}(\Delta)$ for the submanifold of such
$u\in \calo^{1, \alpha}_{\omega ,J}(\Delta)$ that $\re
u|_{[-1.1]}\equiv 0$. Its diffeomorphic image under $\Phi_J$ we denote as 
$\calm_{J,0}\deff\Phi_J\left(\calo^{1,\alpha}_{\omega ,J,0}(\Delta)\right)$.
By $R_{\d\Delta^+}$ denote ``taking real part
and restriction to $\d\Delta^+$'' operator. One has the following commutative diagram:
\begin{equation}
\def\normalbaselines{\baselineskip20pt\lineskip3pt \lineskiplimit3pt }
\def\mapright#1{\smash{\mathop{\longrightarrow}\limits^{#1}}}
\def\mapdown#1{\Big\downarrow\rlap{$\vcenter{\hbox{$\scriptstyle#1$}}$}}
\begin{array}{ccccccccc}
\calc^{1,\alpha}_\omega(\Delta)&\supset&\calo^{1,
\alpha}(\Delta)&\mapright{\Phi_J^{-1}}&\calo^{1,
\alpha}_{\omega ,J}(\Delta)& \\
& &\uparrow{i}& &\uparrow{i}& \\
\calc^{1,\alpha}_\omega(\Delta)&\supset&\calm_{J,0}&\mapright{\Phi_J^{-1}}&
\calo^{1, \alpha}_{\omega
,J,0}(\Delta)&\mapright{R_{\d\Delta^+}}&\calc^{1,
\alpha}_0({\d\Delta^+})&, \\
\end{array}
\eqqno(10)
\end{equation}
where both $i$-s are natural imbeddings. For an unknown map $v$ from $\calo^{1,\alpha}_{\omega
,J,0}(\Delta)$ and given $\phi \in {\calc}^{1,
\alpha}_0({\d\Delta^+})$ consider the system

\begin{eqnarray}
\label{system3}
\begin{cases}
R_{\d\Delta^+}v = \phi , \\
\im v(0) = a.
\end{cases}
\end{eqnarray}

Fix a $J\st$-holomorphic map $v^0\in \calo^{1,\alpha}(\Delta)$ such
that $\re v^0|_{[-1,1]}\equiv 0$ and set $\phi^0=\re
v^0|_{\d\Delta^+}$.
\begin{lem}
\label{analyticity2} For real analytic $J$ close enough to $J\st$ in
$\calc^{1, \alpha}$-norm there exists a neighborhood $U$ of $\phi^0$
in ${\calc}^{1, \alpha}_0({\d\Delta^+})$, a neighborhood ${ U}'$ of
$a^0:=\im v^0(0)$ in $\rr^n$  and a neighborhood $V$ of $v^0$ in
$\calc^{1, \alpha}_{\omega}(\Delta)$ such that for $\phi\in U$ and $a\in U'$
the system (\ref{system3}) admits a unique solution $v \in V\cap
\calo^{1,\alpha}_{\omega ,J,0}(\Delta)$.
\end{lem}

\proof Setting $v=\Phi^{-1}_Jw$ we replace (\ref{system3}) by the  boundary value problem for
an unknown map $w\in \calm_{J,0}$, \ie we shall write it as 

\begin{eqnarray}
\label{system4}
\begin{cases}
R_{\d\Delta^+}\Phi_J^{-1}w = \phi , \\
\im \Phi_J^{-1}w(0) = a.
\end{cases}
\end{eqnarray}

At $J=J\st$ that $\Phi_J=\id$ and $\calm_{J\st , 0}=\calo_0^{1,\alpha}(\Delta)\deff \{ w\in \calo^{1,\alpha}(\Delta):
\re w|{[-1,1]}\equiv 0\}$. The surjectivity condition for the operator obtained by the
linearization of (\ref{system4}) at $w^0 = v^0$ and $J=J\st$ and of the operator
$\Phi^{-1}_J$ is reduced to the resolution of the following system
\begin{equation}
\label{system5} 
\begin{cases}
R_{\d\Delta^+}\dot w  = \psi , \\
\im\dot w (0) = a,
\end{cases}
\end{equation}
for an arbitrary given function $\psi \in {\calc}^{1,
\alpha}_0(\d\Delta^+)$, arbitrary  $a\in \rr^n$ and an unknown map
$\dot w \in \calo^{1, \alpha}_0(\Delta)$. By Lemma \ref{reflectionlemma}
we obtain a solution for any given right hand side of
(\ref{system5}). The uniqueness of $\dot w$ is obvious. Therefore the linearization of (\ref{system4}) 
is a bijective operator at $J=J\st$. By continuity it will be bijective from $T_0\calm_{J,0}$ to 
$\calc_0^{1,\alpha}(\d\Delta^+)\oplus \rr^n$ also for $J$ close to $J\st$. Now the
Implicit Function Theorem implies the desired statement.

\qed

\newprg[prg5.2]{General case} Now we prove Theorem 1.  According to Theorem 2 (which we will prove in the next section) the map $u$ is ${\calc}^\infty$ smooth up to the arc $\beta$. We replace the unit disc by the
upper semi-disc ${\Delta}^+$ and $\beta$ by the interval
$(-1,1)$. By the classical results on the interior regularity of
pseudo-holomorphic maps we can assume that the map $u$ is real
analytic in  a neighborhood of ${ \bar\Delta}^+ \backslash [-1,1]$. Furthermore, the statement is local 
so shrinking $\Delta^+$ if necessary we  assume that $u$ is of class ${\calc}^\infty(\overline\Delta^+)$.
 We can also assume that $X$ is the unit ball of $\cc^n$ equipped with a real
 analytic almost complex structure $J$ with $J(0) = J_{st}$ and that $W  =
 i\rr^n$ and $u(0) = 0$.

Since our considerations are local, we can reduce them to the case of a small deformation of the standard complex structure. 
Indeed, our map $u$ admits
the expansion $u(\zeta) = b\zeta + o(\vert \zeta \vert)$ near the
origin.  For $t > 0$ consider the real analytic structures $J_t(z) =
J(tz)$. They tend to $J_{st}$ as $t$ tends to $0$. Maps $ u^t(\zeta)
= t^{-1} u(t\zeta)$ are $J_t$-holomorphic and tend to the map $
u^0: \zeta \longrightarrow b\zeta$ as $t \longrightarrow 0$ which is
viewed as a  $J_{st}$-holomorphic map.

For $t$ small enough  let $v^t$ be the solution of the boundary-value problem (\ref{system3}) with $\phi^t:= \re u^t \vert
\partial\Delta^+$ given by  Lemma \ref{analyticity2}. This solution is
unique in the class of solutions real analytically extendable past $(-1,1)$. 
However, this still does not give the desired analyticity of $u^t$
since this boundary-value problem could admit solutions in other smoothness classes of maps. So we need the 
uniqueness statement of Lemma \ref{RH-smooth}. Let $H: \Delta \to \Delta^+$ be a biholomorphic map fixing 
the points $-1$ and $1$. Then $H$ is of class ${\calc}^{1/2}(\overline\Delta)$ and extends analytically through 
the open upper and lower semi-circles.
In particular, $H \in L^{1,p}(\Delta)$ for any $p < 4$.  Since $v^t \in {\calc}^{1,\alpha}(\overline\Delta^+)$  
the composition $v^t \circ H$ is in $L^{1,p}(\Delta)$. At the same time $\phi^t \circ H$ belongs to $T^{1,p}$ for every $p < 4$. Applying the 
uniqueness statement of Lemma \ref{RH-smooth} we conclude that $v^t \circ H = u^t \circ H$ which  implies 
that the maps $ u^t$ are real analytic up to $(-1,1)$ for $t$ small enough. This
proves finally that $u$ extends as a real analytic map past
$(-1,1)$. Since it satisfies the real analytic condition
\eqqref(J-holomorphy1) on an open set, the extension is a
$J$-holomorphic map. This proves the Theorem 1.

\section[6]{Compactness}

\newprg[prg6.1]{Compactness theorem: definitions}

We start with recalling some notions and definitions relevant to the
formulation of Gromov compactness theorem as it is stated in
\cite{IS1,IS2}. For a more detailed exposition we refer to these
papers.

\medskip\noindent{\slsf 1. Structures.} We fix a Riemannian manifold
$(X,h)$, a compact subset $K\comp X$ and a sequence of
almost-complex structures $J_n$ of class $\calc^{k,\alpha}$ on $X$,
which converge on $K$ in $\calc^{k,\alpha}$-topology to an
almost-complex structure $J$, $k\ge 0, 0<\alpha <1$. The latter is
supposed to be defined on the whole $X$. Areas of $J_n$-complex
curves will be measured with respect to the "hermitizations"
$h_{J_n}(\cdot ,\cdot ):=\frac{1}{2}(h(\cdot ,\cdot ) + h(J_n\cdot
,J_n\cdot ))$ of $h$ (which converge to $h_J$).

\medskip\noindent{\slsf 2. Curves.} We are given a sequence $\{ C_n\} $
of {\slsf nodal curves with boundary}, parameterized by the same
compact, oriented real surface $(\Sigma ,\d \Sigma )$ with boundary.
We suppose that some parameterizations $\delta_n:\bar\Sigma \to \bar
C_n$ are given as well as some $J_n$-holomorphic maps $u_n:C_n\to X$
such that $(C_n,u_n)$ becomes a {\slsf stable curve over} $(X,J_n)$,
see  Definitions 2.1 - 2.3 from \cite{IS2}.

\medskip\noindent{\slsf 3. Boundedness of areas.} We suppose that
$h_{J_n}$-areas of $u_n(C_n)$ are uniformly bounded and that
$u_n(C_n)\subset K$ for all $n$.

\medskip\noindent{\slsf 4. Weak transversality.} Let $f: W \to X$ be an
immersion of class $\calc^{k+1,\alpha}$ of a real $n$-dimensional
manifold $W$ into a real $2n$-dimensional manifold $X$ and $x\in
f(W)$ a point of self-intersection, so that $f^{-1}(x)= \{w_1, ...,
w_d\} \subset W$ with $d\ge2$. We say that $f(W)$ has {\slsf weakly
transverse self-intersection in $x$} if there exist neighborhoods
$U_i \subset W$ of $w_i$ such that for any pair $w_i \not= w_j$ the
intersection $f(U_i) \cap f(U_j)$ is a
$\calc^{k+1,\alpha}$-submanifold in $X$ of dimension equal to $\dim
\bigl( df(T_{w_i}W) \cap df(T_{w_j}W) \bigr)$.

\medskip\noindent{\slsf 5. Totally real submanifolds.} We fix a
collection $\{ W_i\}_{i=1}^m$ of real manifolds of real dimension
$n= \dimc X$ and $J$-totally real immersions  $f_i:W_i\to X$, i=1,...,m,
of class $\calc^{k+1,\alpha}$.  We call ${\mib W}=
\{(W_i, f_{i})\}_{i=1}^m$ the $J$-totally real immersed submanifold
of $X$. We suppose that the immersion ${\mib f} =\{f_i\}$ has only
{\slsf weakly transverse} self-intersections. More precisely that
means that each $f_i$ has only weakly transverse self-intersections
and every pair $f_i, f_j$ intersect in the same manner.

Furthermore, we fix a sequence ${\mib W}_n = \{(W_i,
f_{n,i})\}_{i=1}^m$ of $J_n$-totally real immersed submanifolds
(with weakly transverse self-intersections), which  converge in
$\calc^{k+1,\alpha}$-sense to a $J$-totally real submanifold ${\mib
W} =\{ ( W_i,f_{i})_{i=1}^m\} $ (again, with weakly transverse
self-intersections!).

\begin{rema}\rm The condition of a weak transverse self-intersection
is crucial in applications. The reason is that if, for example, a
totally real manifold $W$ is immersed into $\cc^n$ with transversal
self-intersections then its product with, say a circle $W\times
\ss^1$ will be immersed into $\cc^{n+1}$, but its self-intersections
will be only weakly transverse! This construction repeatedly occurs
in applications. See more about this in \cite{IS2}.
\end{rema}

\medskip\noindent{\slsf 6. Boundary conditions.} We fix a collection of
arcs with disjoint interiors $\bfbeta= \{\beta_k\}_{k=1}^M$, which
defines a decomposition of the boundary $\d \Sigma = \cup_k
\beta_k$. We assume that every boundary point $b$ of $\Sigma $,
which is mapped by the parametrization $\delta_n :\bar\Sigma \to
C_n$ into a boundary node $a$ of $C_n$, is the endpoint of two arcs
and that $a$ itself is an endpoint for four arcs $\delta_n (\beta_k
)=:\beta_{n,k}$. Our basic assumption is that the same collection
$\bfbeta $ serves for {\slsf all} curves $C_n$. Totally real boundary
conditions $({\mib W}, \bfbeta ,{\mib u}^{(b)})$ is the data, which
includes ${\mib W}=\{(W_i, f_i)\}$, $\bfbeta =\{\beta_k\}$ and
continuous maps ${\mib u}^{(b)}  =\{ u^{(b)}_{k} : \beta_k \to W_i\}
, i=1,...,m, k=1,...,M$. Here several different $\beta_k$-s could be
mapped into the same $W_i$.

We shall suppose that the sequence of curves $(\barr C_n, u_n)$ satisfy the {\slsf
totally real boundary conditions} $({\mib W}_n,\bfbeta ,{\mib
u}^{(b)}_n )$ in the sense  that there are given continuous maps
$u^{(b)}_{n,k} :\beta_k \to W_i$ with $f_{n,i} \circ u^{(b)}_{n,k} =
u_n|_{\beta_{n,k}}=u_n\circ \delta_n\mid_{ \beta_k}$. Here ${\mib
u}^{(b)}_n =\{ u^{(b)}_{n,k}\} $.

\medskip\noindent{\slsf 7. Description of convergency.} Compactness
theorem states that under the assumptions described above there
exist a subsequence $(C_{n_k},u_{n_k})$ which converge in the
following sense.

\begin{defi}
\label{converge} We say that the sequence $(\barr C_n, u_n)$ of
stable $J_n$-complex curves over $X$, which satisfies the totally
real boundary conditions $({\mib W}_n,\beta, {\mib u}_n^{(b)})$ {\slsf
converges up to the boundary} to a stable $J$-complex curve $(\barr
C, u)$ over $X$ if the parameterizations $\sigma_n: \barr\Sigma \to
\barr C_n$ and $\sigma : \barr\Sigma \to \barr C$ can be chosen in
such a way that the following holds:

\sli $u_n\scirc \sigma_n$ converges to $u\scirc \sigma$ in $C^0(
\barr\Sigma, X)$-topology; moreover, $u_n^{(b)}$ uniformly on $\d
\Sigma $ converge to some $u^{(b)}$ such that $(\bar C,u)$ satisfies
the totally real boundary condition $({\mib W},\bfbeta, {\mib
u}^{(b)} )$;

\slii if $\{ a_k \}$ is the set of the nodes of $C_\infty$ and
$\gamma_k \deff \sigma\inv(a_k)$ are the corresponding circles and
arcs in $\barr\Sigma$, then for any compact subset $R\comp
\barr\Sigma \bs \cup_k\gamma_k$ there exists $n_0=n_0(K)$, such that
$ \sigma_n^{-1}(\{ a_k \}) \cap K= \emptyset$ for all $n\ge n_0$ and
complex structures $\sigma_n^*j_{C_n}$ smoothly converge to
$\sigma^*j_{C}$ on $R$, $n\ge n_0$;

\sliii  on any compact subset $R\comp \barr\Sigma \bs
\cup_k\gamma_k$ the convergence $u_n\scirc \sigma_n\to u\scirc
\sigma$ is in $\calc^{k+1,\alpha}$-topology.
\end{defi}

\smallskip
Note that in this definition it is {\sl ad hoc} supposed that all
$\bar C_n$-s and also $\bar C$ can be parameterized by the same real
surface $\bar \Sigma $ and that one is allowed to change
parameterizations from initial $\delta_n$ to an appropriate
$\sigma_n$. Note also that no convergence of ${\mib u}_n^{(b)} )$ is
a priori supposed. It comes as the statement of the Theorem.

\newprg[6.2]{Generalized Giraud-Calderon-Zygmund Inequality}

 Recall the following Giraud Inequality (or estimate):
 \it for all $1<p<\infty $ there exists a constant $G_p$
 such that for all $u\in L^p(\Delta ,\cc^n )$ one has

\begin{equation}
\norm{ (\d \scirc T_{\Delta}^{CG})(u)}_{L^p(\Delta)}\leq G_p\cdot
\norm{u}_{L^p(\Delta)}. \label{G-est}
\end{equation}
\rm In H\"older norms an analogous statement is due to Calderon and
Zygmund. Namely: \it for every $0<\alpha <1$ there exists
$C_{\alpha}$ such that for all $u\in \calc^{\alpha}(\Delta , \cc^n)$
one has
\begin{equation}
\norm{ (\d \scirc T_{\Delta}^{CG})(u)}_{L^{\alpha}(\Delta)}\leq
G_{\alpha}\cdot \norm{u}_{L^{\alpha}(\Delta)}. \label{CZ-est}
\end{equation}
\rm We shall need the following generalization of (\ref{G-est}) and
(\ref{CZ-est}) to $\dbar$-type operators. Let $J=J(\zeta)$ be a
bounded (resp. $\calc^{\alpha}$-continuous) operator in the trivial
bundle $\Delta\times\rr^{2n}$ which satisfies $J^2(\zeta)\equiv
-\id$ for all $\zeta\in\Delta$. It defines a natural $\dbar$-type
operator $\dbar_J: L^{1,p}(\Delta,\rr^{2n})\to L^p(\Delta,\rr^{2n})$
(resp. $\dbar_J: \calc^{1,\alpha}(\Delta,\rr^{2n})\to
\calc^{\alpha}(\Delta,\rr^{2n})$) as follows
\begin{equation}
\eqqno(debar-J) \dbar_Ju = \frac{\d u}{\d\xi} + J(\zeta)\frac{\d
u}{\d\eta}.
\end{equation}

\begin{lem}
\label{GCZ-gen} For any $p>2$ there exist $\eps_p>0$ and
$C(p,\norm{J-J\st}_{L^{\infty}})<\infty$ $\big($resp. for any
$0<\alpha <1$ there exist $\eps_{\alpha}>0$ and $C(\alpha,
\norm{J-J\st}_{\alpha})\big)$ such that for any $J\in L^{\infty}(\Delta , End (\rr^{2n}))$,
$J^2\equiv -\id$ with
$\norm{ J - J\st }_{L^{\infty}(\Delta)} <\eps_p$ $\big($resp. any
$J\in\calc^{\alpha}$ with
$\norm{J-J\st}_{\calc^{\alpha}}<\eps_{\alpha}\big)$ any $u\in
L^{p}(\Delta,\rr^{2n})$ $\big($resp. any $u\in
\calc^{\alpha}(\Delta,\rr^{2n})\big)$ with compact support in
$\Delta$ one has

\begin{equation}
\eqqno(G-gen) \norm{du}_{L^p(\Delta ,\rr^{2n})}\le
C(p,\norm{J-J\st}_{L^{\infty}})\norm{\dbar_Ju}_{L^p(\Delta
,\rr^{2n})},
\end{equation}
\noindent and respectively
\begin{equation}
\eqqno(CZ-gen) \norm{du}_{\calc^{\alpha}(\Delta ,\rr^{2n})}\le
C(\alpha ,\norm{J-J\st}_{\calc^{\alpha}})
\norm{\dbar_Ju}_{\calc^{\alpha}(\Delta ,\rr^{2n})}.
\end{equation}
\end{lem}
\proof For the proof of \eqqref(G-gen) see Lemma 1.2 in \cite{IS1}.
The proof of \eqqref(CZ-gen) follows the same lines.  For $u\in
\calc^{\alpha}(\cc ,\rr^{2n})$ it holds that
\[
\norm{ (\dbar_J\scirc T_{\Delta}^{CG} - \dbar_{J\st}\scirc
T_{\Delta}^{CG})u}_{\calc^{\alpha}(\Delta )}\leq \norm{ J -
J\st}_{\calc^{\alpha}(\Delta )}\cdot \norm{
d(T_{\Delta}^{CG}u)}_{\calc^{\alpha}(\Delta )} \leq
\]
\begin{equation}
\leq \norm{ J - J\st}_{\calc^{\alpha}(\Delta )}(1+G_{\alpha})\norm{
u}_{\calc^{\alpha}(\Delta )}, \label{c-alfa-ozenka}
\end{equation}
where $G_{\alpha}$ is the constant from (\ref{CZ-est}). For  the
standard structure in $\cc^n$ the operator  $\dbar_{J\st}\scirc
T_{\Delta}^{CG}: \calc^{\alpha}(\Delta ,\cc^n)\to
\calc^{\alpha}(\Delta,\cc^n)$ is the identity. So from
(\ref{c-alfa-ozenka}) we see that there exists
$\eps_{\alpha}=\frac{1}{1+G_{\alpha}}$ such that if $\norm{ J - J\st
}_{\calc^{\alpha}(\Delta)} <\eps_{\alpha}$, then $\dbar_{J}\scirc
T_{\Delta}^{CG}: \calc^{\alpha}(\Delta, \cc^n) \to
\calc^{\alpha}(\Delta , \cc^n)$ is an isomorphism. Moreover, since
$\dbar_J\scirc T_{\Delta}^{CG} = \dbar_{J\st}\scirc T_{\Delta}^{CG}
+ (\dbar_J-\dbar_{J\st})\scirc T_{\Delta}^{CG}$, we have
\begin{equation}
(\dbar_J\scirc T_{\Delta}^{CG})\inv = \left[ \id +
(\dbar_J-\dbar_{J\st})\scirc T_{\Delta}^{CG}\right] ^{-1} =
\sum_{n=0}^{\infty }(-1)^n[(\dbar_J - \dbar_{J\st})\scirc
T_{\Delta}^{CG}]^n.
\end{equation}
This shows, in particular, that $(\dbar_J\scirc T_{\cc}^{CG})\inv $
does not depend on $0<\alpha< 1$ provided that $\norm{ J - J\st
}_{\calc^{\alpha}(\Delta)} <\eps_{\alpha}$.

\smallskip Put $h = u - T_{\Delta}^{CG}\scirc \dbar_{J\st}u$. Then
$\dbar_{J\st}h=0$. So $h$ is holomorphic and descends at infinity.
Thus $h\equiv 0$, which implies $u=(T_{CG}\scirc \dbar_{J\st})u$.
Write
\[
\norm{du}_{\calc^{\alpha}(\Delta)}\le (1+C_{\alpha})\norm{
\dbar_{J\st }u}_{\calc^{\alpha}(\Delta )} =
(1+C_{\alpha})\norm{(\dbar_J\circ T^{CG}_{\Delta})^{-1}(\dbar_J\circ
T^{CG}_{\Delta}) \dbar_{J\st }u}_{\calc^{\alpha}(\Delta )} =
\]
\[
= (1+C_{\alpha})\norm{(\dbar_J\circ
T^{CG}_{\Delta})^{-1}(\dbar_Ju}_{\calc^{\alpha}(\Delta )} \leq
(1+C_{\alpha})\sum_{n=0}^{\infty }\norm{
(\dbar_J-\dbar_{J\st})\scirc
T^{CG}_{\Delta}}_{\calc^{\alpha}(\Delta)}^n \cdot \norm{
\dbar_Ju}_{\calc^{\alpha}(\Delta )} \leq
\]
\begin{equation}
\leq  C(\alpha , \norm{J-J\st}_{\calc^{\alpha}})\cdot \norm{
\dbar_Ju}_{\calc^{\alpha}(\Delta )},
\end{equation}
provided that $\norm{J-J\st}_{\calc^{\alpha}}<\eps_{\alpha} $.

\smallskip\qed

\begin{corol}
\label{c1-alfa-conv}
If $J_n\to J$ in $\calc^{k,\alpha}$-norm on compact $K\comp X$,
$k\ge 0, 0<\alpha <1$, and $J_n$-holomorphic maps $u_n:\Delta\to
K\comp X$ uniformly converge to $u:\Delta\to X$ then $u_n$
converge to $u$ in $\calc^{k+1,\alpha}$-topology on compacts in $\Delta$.
\end{corol}
\proof This will be done in three steps.

\smallskip\noindent{\slsf Step1. $\calc^{\alpha}$-convergence.}
First we prove the $\calc^{\alpha}$-convergency (which is
implicitly contained in \cite{IS1}). For this we need only {\slsf
uniform} convergence of {\slsf continuous} structures. Consider all
$u_n$ as a sections of the trivial bundle $\Delta\times\rr^{2n}$
which are holomorphic with respect to the pulled back structures
$J_n\circ u_n$, \ie $\dbar_{J_n\circ u_n}u_n=0$. Theorem 6.2.5 from
\cite{M} implies that for  every $2\le p<\infty$

\begin{equation}
\eqqno(morrey) \norm{u_n}_{L^{1,p}(\Delta (1/2)} \le C_p\big(\mu
(J_n\circ u_n)\big)\norm{u_n}_{L^2(\Delta)},
\end{equation}
where the constant $C_p\big(\mu (J_n\circ u_n)\big)$ crucially
depends not only on $p$ but also on the modulus of continuity of
$J_n\circ u_n$. This gives us the boundedness of $u_n$ in
$L^{1,p}(\Delta (1/2))$ for all $p$ and therefore in
$\calc^{\gamma}$ for all $\gamma <1$ by Sobolev imbedding theorem.
And the last in its turn by the Ascoli theorem implies the
$\calc^{\gamma}$-convergency for all $\gamma$, in particular, for
our $\alpha$ in question.

\smallskip\noindent{\slsf Step 2. $\calc^{1,\alpha}$-convergence.} 
Take a cut-off function $\phi$
and write using \eqqref(CZ-gen) (all constants below are denoted by the same letter
$C$, but are different):
\[
\norm{d\big[\phi (u_n - u_m)\big]}_{\calc^{\alpha}}\le
C\norm{\dbar_{J_n\circ u_n}\big[\phi (u_n -
u_m)\big]}_{\calc^{\alpha}}\le C\norm{u_n-u_m}_{\calc^{\alpha}} +
C\norm{\phi\dbar_{J_n}u_m}_{\calc^{\alpha}} \le
\]
\[
\le C\norm{u_n-u_m}_{\calc^{\alpha}} + \norm{J_n\circ u_n-J_m\circ u_m}_{\calc^{\alpha}}\norm{\phi
du_m}_{\calc^{\alpha}}\le C\norm{u_n-u_m}_{\calc^{\alpha}} + 
\]
\begin{equation}
\eqqno(5.10) 
+ C\norm{J_n\circ u_n-J_m\circ
u_m}_{\calc^{\alpha}}\big[\norm{d(\phi u_m)}_{\calc^{\alpha}} 
+ \norm{u_m}_{\calc^{\alpha}}\big].
\end{equation}
In the same manner from \eqqref(CZ-gen) we get that
\[
\norm{d(\phi u_m)}_{\calc^{\alpha}}\le C\norm{\dbar_{J_m\circ
u_m}(\phi u_m)}_{\calc^{\alpha}} \le C\norm{u_m}_{\calc^{\alpha}},
\]
and the latter is bounded. \eqqref(5.10) implies now
$\calc^{1,\alpha}$-convergence of $u_n$ to $u$.

\smallskip\noindent{\slsf Step 3. $\calc^{k+1.\alpha}$-convergence.} 
We already lifted an
almost complex structure from the manifold $X$ to its tangent bundle $TX$, see \eqqref(lift).
Since $u_n$ already converge in $\calc^{1,\alpha}$-topology by Step 2,
lifts $u_n^c$ converge in $\calc^{\alpha}$-topology. And therefore again by Step 2 they
converge in $\calc^{1,\alpha}$-topology. The rest is obvious. 

\smallskip\qed

\newprg[prg6.3]{Proof of the Compactness theorem}

In \cite{IS2} it was proved that there exists a subsequence
$u_{n_k}$ which converges as in Definition \ref{converge} with the
only  difference that in (\sliii the convergence was in
$L^{1,p}$-topology for all $p<\infty$. That implies
$\calc^{\alpha}$-convergence for all $0<\alpha <1$. We recall that
structures $J_n$ in \cite{IS2} where supposed to be only continuous
and uniformly, \ie $\calc^0$-convergent to $J$. All we need to do
here is to improve convergence to $\calc^{k+1,\alpha}$ in the case
when $J_n$ converge in $\calc^{k,\alpha}$ and $W_n$ in
$\calc^{k+1,\alpha}$. 

\smallskip The statement we need is purely local and therefore the
totally real manifolds can be supposed to be imbedded. $\beta$, as
before, stays for the interval $(-1,1)$.

\begin{lem}
\label{ck-alfa-reg} Let $k\ge 0$ and let a sequence $\{J_n\}$ of almost complex
structures of class $\calc^{k,\alpha}$ converge in
$\calc^{k,\alpha}$-topology to $J$, a sequence $\{W_n\}$ of imbedded
$J_n$-totally real submanifolds of class $\calc^{k+1,\alpha}$
converge in $\calc^{k+1,\alpha}$- topology to a $J$-totally real
$W$, and a sequence $u_n:(\Delta^+,\beta)\to (X,W_n)$ of
$J_n$-holomorphic maps of class $\calc^0\cap L^{1,2}$ converge in
$\calc^0\cap L^{1,2}$-topology to $u:(\Delta^+,\beta)\to (X,W)$.
Then $u_n$ converge to $u$ in $\calc^{k+1,\alpha}$-sense up to
$\beta$.
\end{lem}
\proof We start with the case  $k=0$. An obvious modification of Lemma \ref{alpha}
provides us a sequence $\{\phi^n\}$ of $\calc^{1,\alpha}$-
diffeomorphisms converging in $\calc^{1,\alpha}$-topology to $\phi$
and such that they all satisfy the properties 1)-3) stated there.
Repeating the considerations of the proof of Lemma \ref{c1-alfa-reg} to each $J_n$ and $u_n$ we
get a uniform bound on the norms of ``extensions by reflection'' $\norm{\tilde u_n}_{\calc^{1,\alpha}}$ and
this implies the statement of the Lemma for $k=0$ by Corollary \ref{c1-alfa-conv}.

Having  convergency in $\calc^{1,\alpha}$ topology we obtain the convergency in higher regularity classes via
geometric bootstrap of  subsection 3.3. Really, if $J_n\to J$ in $\calc^{1,\alpha}$-topology  then 
$J_n^c\to J^c$ in $\calc^{\alpha}$-topology. Further we  lift 
$J_n$-holomorphic maps $u_n:\Delta\to X$  to  $J^c$-holomorphic maps
$u^c_n:\Delta\to TX$ and apply the case $k=0$. 

\smallskip\qed

\section[8]{Open Questions}
At the end we would like to turn the attention of a reader to some
open questions.

\smallskip\noindent
{\bf 1.} Let $(X,J)$ be a real analytic almost complex manifold and
$W$ a real analytic $J$-totally real submanifold of $X$. Let $C^+$
be $J$-complex curve in $X\setminus W$. Does there exists a
neighborhood $V$ of $W$ and a $J$-complex curve $C^-$ in $V\setminus
W$ (reflection of $C^+$) such that $\overline{\left( C^+\cup
C^-\right) }\cap V$ is a $J$-complex curve in $V$?

\medskip For  integrable $J$ the answer is yes and is due to H.
Alexander, see \cite{A}.

\medskip\noindent
{\bf 2.} The following question is a particular case of the previous
one. Let $C$ be a $J$-complex curve in the complement of a point.
Will its closure $\bar C$ be a $J$-complex curve?

\medskip\noindent
{\bf 3.} This question was communicated to us by J.-C. Sikorav.
Define a $J$-holomorphic map as a differentiable map $u:\Delta\to X$
such that \eqqref(J-holomorphy1) is satisfied at every point. Prove
that $u\in L^{1,2}_{loc}$ (and therefore $u$ is a $J$-holomorphic
map in the usual sense).

\ifx\undefined\bysame
\newcommand{\bysame}{\leavevmode\hbox to3em{\hrulefill}\,}
\fi

\def\entry#1#2#3#4\par{\bibitem[#1]{#1}
{\textsc{#2 }}{\sl{#3} }#4\par\vskip2pt}

\bigskip

\ \vspace{-30pt}

\newdimen\widthone
            \newdimen\widthtwo
            \newdimen\widthcommon
            \def\twocol#1#2{
            \widthcommon=\hsize
            \setbox1=\hbox{Hochschulausbilbung mmmmmmmmmmmmmm\ \ \ \hskip 1.5pt}
            \widthone=\wd1\relax
            \widthtwo=\widthcommon \advance\widthtwo by -\widthone
            \advance\widthtwo by -.3 true cm
            \hbox{\hsize=\widthcommon\noindent
            {\vtop{\hsize=\widthone\noindent\sl#1}}%
            {{\vtop{\hsize=\widthtwo\hbadness=1600\noindent#2}}}}\par}

 \twocol{U.F.R. de Math\'ematiques, \\
        Universit\'e de Lille-1\\
        59655 Villeneuve d'Ascq, France.\\
         E-mail: ivachkov@math.univ-lille1.fr\\ \\
         IAPMM Acad. Sci. Ukraine,\\
         Lviv, Naukova 3b,\\
         79601 Ukraine.}
         {U.F.R de Math\'ematiques,\\
         Universit\'e de Lille-1,\\
         59655 Villeneuve d'Ascq, France,\\
         E-mail: sukhov@math.univ-lille1.fr}

\end{document}